%% file: main.tex
\documentclass[11pt,letterpaper]{amsart}

\usepackage{graphicx}
\usepackage{mathtools,amsmath,amsthm,amsfonts,amssymb}
\usepackage{enumitem}
\setlist[enumerate]{font=\normalfont}
\usepackage{tikz}
\usepackage{verbatim}
\usetikzlibrary{patterns}

\usepackage[colorlinks,allcolors=blue,backref=page]{hyperref}
\usepackage[capitalize]{cleveref}

\newtheorem{theorem}{Theorem}[section]
\newtheorem{lemma}[theorem]{Lemma}

\newtheorem{proposition}[theorem]{Proposition}
\newtheorem{corollary}[theorem]{Corollary}

\theoremstyle{definition}
\newtheorem{definition}[theorem]{Definition}

\newtheorem{example}[theorem]{Example}

\theoremstyle{remark}
\newtheorem{remark}[theorem]{Remark}

\renewcommand*{\backrefalt}[4]{%
  \ifcase #1 %
  \or
    ↑~#2%
  \else
    ↑~#2%
  \fi}

\usepackage{xcolor}
\definecolor{airforceblue}{rgb}{0.36, 0.54, 0.66}

\definecolor{zsxpurple}{RGB}{128,0,128}

\definecolor{Lightred}{RGB}{239,83,80}

\title[Kazhdan--Lusztig Polynomials with \(P_{\mathrm{id},w}(1)=3\)]
{Kazhdan--Lusztig Polynomials for permutations with \(P_{\mathrm{id},w}(1)=3\) and Irreducible Singular Locus}

\author[Gaston Burrull]{Gaston Burrull}
\address{International Institute for Mathematical Sciences, Jilin University, Qianjin Avenue 2699, High-Tech Development Zone, Changchun, Jilin 130012, China}
\email{grburrul@gmail.com}

\author[Zhitong Alex Li]{Zhitong Alex Li}
\address{Department of Mathematics, Dornsife College of Letters, Arts and Sciences, University of Southern California, Los Angeles, CA 90089, USA}
\email{ali24951@usc.edu}

\author[Shixuan Zeng]{Shixuan Zeng}
\address{School of the Gifted Young, University of Science and Technology of China, No. 96 Jinzhai Road, Baohe District, Hefei, Anhui 230026, China}
\email{zsx0515@mail.ustc.edu.cn}

\date{}

\begin{document}

\input{0_abstract}
\maketitle
\input{1_introduction}
\input{acknowledgments}
\input{2_preliminaries}
\input{3_proof_R1}
\input{4_proof_R2}
\input{5_proof_R3}

\end{document}

%% file: 0_abstract.tex
\begin{abstract}
    We determine the Kazhdan--Lusztig polynomial \(P_{\mathrm{id},w}(q)\) for permutations \(w\) such that \(P_{\mathrm{id},w}(1)=3\) and the Schubert variety \(X_w\) has irreducible singular locus.
    Let \(h=h(w)\) be the minimum height \(w(a)-w(d)\) of a \(3412\)-occurrence in \(w\) at positions \(a<b<c<d\), with \(h=1\) when no such occurrence exists.
    Then
    \[
        P_{\mathrm{id},w}(q)
        =
        \begin{cases}
            1+q^h+q^{h+1},
            &
            \text{if exactly one of $3412$ and $4231$ occurs},\\[2mm]
            1+q+q^h,
            &
            \text{if both patterns occur}.
        \end{cases}
    \]
    In the second case, \(h>1\).
    This proves the irreducible singular locus case of a 2009 conjecture of Billey and Weed.
\end{abstract}

%% file: 1_introduction.tex
\section{Introduction}\label{sec:introduction}

Kazhdan--Lusztig polynomials, introduced in~\cite{KL79}, play a central role in representation theory and are closely connected to the geometry of Schubert varieties.
For $w\in S_n$, let $X_w$ denote the corresponding Schubert variety in the flag variety of $\mathrm{GL}_n(\mathbb C)$.
For $u\leq w$ in Bruhat order, the Kazhdan--Lusztig polynomial $P_{u,w}(q)$ records the local intersection cohomology of $X_w$ along the Schubert cell $X_u^{\circ}$~\cite{KL80}.
In particular, $P_{\mathrm{id},w}(1)$ gives the total dimension of the local intersection cohomology at the identity fixed point, while $P_{\mathrm{id},w}(q)$ determines its graded dimensions.

Although Kazhdan--Lusztig polynomials can be computed recursively~\cite[Section~5.1]{BB05}, obtaining explicit formulas remains difficult even in type~$A$.
Lascoux and Sch\"utzenberger gave a combinatorial formula for Kazhdan--Lusztig polynomials of cograssmannian permutations~\cite{LS81}, and Lascoux extended this to covexillary permutations~\cite{Las95}.
From a geometric perspective, Billey and Warrington classified the irreducible components of the singular loci of Schubert varieties and determined the Kazhdan--Lusztig polynomials at their generic points~\cite{BW03}.

Another direction is to classify permutations for which $P_{\mathrm{id},w}(1)$ is small and determine the corresponding polynomials $P_{\mathrm{id},w}(q)$.
Since $P_{\mathrm{id},w}(q)$ has nonnegative integer coefficients and constant term one, the smallest possible value of $P_{\mathrm{id},w}(1)$ is $1$.
By Deodhar's criterion~\cite{Deo85}, this occurs precisely when $X_w$ is smooth, or equivalently, by the Lakshmibai--Sandhya criterion~\cite{LS90}, when $w$ avoids both $3412$ and $4231$.

The next case, $P_{\mathrm{id},w}(1)=2$, was studied by Billey and Braden~\cite{BBr03}.
They showed that this condition implies irreducibility of the singular locus of $X_w$ and avoidance of six particular patterns of length six, and conjectured the converse.

For a \(3412\)-embedding in \(w\) at positions \(a<b<c<d\), its \emph{height} is \(w(a)-w(d)\).
Let \(h(w)\) be the minimum height of such embeddings, with \(h(w)=1\) when \(w\) avoids \(3412\).
Woo proved Billey and Braden's conjecture and obtained the formula~\cite[Theorem~1.1]{Woo09}
\[
    P_{\mathrm{id},w}(q)=1+q^{h(w)}.
\]

For \(P_{\mathrm{id},w}(1)=3\), nonnegativity and the constant term condition give \(P_{\mathrm{id},w}(q)=1+2q^r\) for some \(r>0\), or \(P_{\mathrm{id},w}(q)=1+q^r+q^s\) for some \(0<r<s\).
In the appendix to Woo's paper, Billey and Weed conjectured that the second possibility holds whenever the singular locus of $X_w$ is irreducible~\cite[Conjecture~A.5(2)]{Woo09}.
We prove this conjecture by determining \(r\) and \(s\) explicitly.

\begin{theorem}
\label{thm:main}
    Let $w\in S_n$ satisfy $P_{\mathrm{id},w}(1)=3$,
    and suppose that the singular locus of $X_w$ is irreducible.
    Set $h:=h(w)$.
    Then
    \[
        P_{\mathrm{id},w}(q)
        =
        \begin{cases}
            1+q^h+q^{h+1},
            &
            \text{if $w$ contains exactly one of $3412$ and $4231$},\\[2mm]
            1+q+q^h,
            &
            \text{if $w$ contains both $3412$ and $4231$}.
        \end{cases}
    \]
    Moreover, $h\geq2$ in the second case.
\end{theorem}

For singular \(X_w\), let \(d(w)\) denote the smallest positive exponent of \(P_{\mathrm{id},w}(q)\).
Gaetz and Gao determined this exponent as follows.

\begin{theorem}[{\cite[Theorem~1.6]{GG24}}]
\label{thm:Gaetz-Gao-min-degree}
    Let $w\in S_n$ satisfy $P_{\mathrm{id},w}(q)\neq1$.
    Then
    \[
        d(w)=
        \begin{cases}
            1, & \text{if $w$ contains $4231$},\\
            h(w), & \text{if $w$ avoids $4231$}.
        \end{cases}
    \]
\end{theorem}

In the notation above, Gaetz--Gao's theorem determines \(r\), while Theorem~\ref{thm:main} also determines \(s\) under stronger hypotheses.

\subsection*{Overview of the proof}
Under the hypotheses of Theorem~\ref{thm:main}, \(X_w\) is singular and therefore \(w\) contains at least one of \(3412\) and \(4231\).
Set \(h:=h(w)\).
The proof splits into three cases:

\begin{enumerate}[
        label=\textnormal{(R\arabic*)},
        ref=\textnormal{(R\arabic*)},
        leftmargin=*
    ]
    \item\label{case:R1}
    If \(w\) contains \(3412\) and avoids \(4231\), then \(P_{\mathrm{id},w}(q)=1+q^h+q^{h+1}\).

    \item\label{case:R2}
    If \(w\) avoids \(3412\) and contains \(4231\), then \(P_{\mathrm{id},w}(q)=1+q+q^2\).

    \item\label{case:R3}
    If \(w\) contains both patterns, then \(h\geq2\) and \(P_{\mathrm{id},w}(q)=1+q+q^h\).
\end{enumerate}


We prove these statements in Sections~\ref{sec:3412-only}, \ref{sec:4231-only}, and~\ref{sec:mixed}, respectively.
For~\ref{case:R1}, we compute the identity fiber of Cortez's resolution (Section~\ref{sec:3412-only}) and apply the decomposition theorem.
For~\ref{case:R2}, we use Lascoux's formula (Theorem~\ref{thm:R2-Lascoux-tree}) and the irreducibility of the singular locus to reduce the calculation to labelings of a single chain.
For~\ref{case:R3}, we first show that $h\geq2$ and use continuous restrictions to reduce to a simpler permutation.
The Billey--Braden inequality (Theorem~\ref{thm:Billey-Braden-product}) then gives a positive coefficient of $q^h$, while Gaetz--Gao's theorem (Theorem~\ref{thm:Gaetz-Gao-min-degree}) gives a positive coefficient of $q$.

%% file: acknowledgments.tex
\subsection*{Acknowledgments}
This project is part of the PACE 2026 program at BICMR, Peking University, whose organizers we thank for their hospitality.
We also thank Yibo Gao and Xinyuan Zhang for useful discussions.

%% file: 2_preliminaries.tex
\section{Preliminaries}\label{sec:preliminaries}

\subsection{Bruhat order and pattern containment}
\label{subsec:bruhat-patterns}

Let $S_n$ be the symmetric group on $[n]:=\{1,\ldots,n\}$, with identity $\mathrm{id}$ and simple reflections
\[
    s_i=(i,i+1),
    \qquad
    1\leq i\leq n-1.
\]
We write permutations in one-line notation \(w=w(1)\cdots w(n)\) and compose them as functions, so that \((uv)(i)=u(v(i))\).
The \emph{length} of $w$ is its number of inversion pairs,
\[
    \ell(w)
    =
    \#\{(i,j):1\leq i<j\leq n,\ w(i)>w(j)\}.
\]

We consider the strong Bruhat order $\leq$ on $S_n$.
For $u\leq w$, we write
\[
    [u,w]:=\{v\in S_n:u\leq v\leq w\}
\]
for the corresponding Bruhat interval.

We will use the classical tableau criterion for Bruhat order.
Its geometric form goes back to Ehresmann \cite[\S\S19--20, pp.~438--442]{Ehr34}; see \cite[p.~1]{BB96} for the permutation formulation used here.

\begin{lemma}
\label{lem:bruhat-tableau-criterion}
    For $w\in S_n$ and $1\leq p\leq n-1$, let
    \[
        w^{(p)}_1<\cdots<w^{(p)}_p
    \]
    be the increasing rearrangement of $w(1),\ldots,w(p)$.
    Then, for
    $u,v\in S_n$,
    \[
        u\leq v
        \quad\Longleftrightarrow\quad
        u^{(p)}_j\leq v^{(p)}_j
        \quad
        \text{for all }1\leq j\leq p\leq n-1.
    \]
\end{lemma}

\begin{definition}\label{def:flattening-pattern-containment-embedding}
    For a sequence $a_1,\ldots,a_k$ of distinct integers, its
    \emph{flattening}
    \[
        \operatorname{fl}(a_1\cdots a_k)\in S_k
    \]
    is the unique permutation with the same relative order as the sequence.
    Given $v\in S_k$ and $w\in S_n$, an \emph{embedding} of $v$ in $w$ (or \emph{$v$-embedding in~$w$}) is a
    sequence of indices
    \(1\leq i_1<\cdots<i_k\leq n\)
    such that
    \[
        \operatorname{fl}\bigl(w(i_1)\cdots w(i_k)\bigr)=v.
    \]
    We say that $w$ \emph{contains} $v$ if such an embedding exists, and that
    $w$ \emph{avoids} $v$ otherwise.
\end{definition}

We will also use the interval version of pattern embeddings.

\begin{definition}
\label{def:interval-pattern-embedding}
    Let $[x,v]\subseteq S_m$ and $[u,w]\subseteq S_n$ be Bruhat intervals.
    An \emph{interval pattern embedding} of $[x,v]$ into $[u,w]$ is a
    simultaneous embedding of $x$ in $u$ and of $v$ in $w$, using the same set of positions \(J\subseteq[n]\), such that
    \begin{enumerate}[label=(\roman*), ref=(\roman*)]
        \item\label{item:interval-pattern-outside}
        \(u(i)=w(i)\) for all \(i\notin J\).
        
        \item\label{item:interval-pattern-isomorphism}
        $[x,v]$ and $[u,w]$ are isomorphic as posets.
    \end{enumerate}
\end{definition}

By \cite[Lemma~2.1]{WY08}, for a simultaneous embedding as above satisfying~\ref{item:interval-pattern-outside}, condition~\ref{item:interval-pattern-isomorphism} is equivalent to
\[
    \ell(v)-\ell(x)=\ell(w)-\ell(u).
\]

\begin{example}
    Consider \(u=14325\) and \(w=52341\) in \(S_5\).
    The tableau criterion gives \(u\leq w\).
    For instance, for \(p=3\),
    \[
        \left(u^{(3)}_1,u^{(3)}_2,u^{(3)}_3\right)=(1,3,4),
        \qquad
        \left(w^{(3)}_1,w^{(3)}_2,w^{(3)}_3\right)=(2,3,5),
    \]
    so \(u^{(3)}_j\leq w^{(3)}_j\) for \(1\leq j\leq3\).

    Let \(x=1324\), \(v=4231\), and \(J=\{1,2,4,5\}\).
    The positions in \(J\) give a simultaneous embedding of \(x\) in \(u\) and \(v\) in \(w\) since \(\operatorname{fl}(1425)=1324\) and \(\operatorname{fl}(5241)=4231\).
    Moreover, \(u(3)=w(3)=3\).
    Finally,
    \[
        \ell(w)-\ell(u)
        =
        4
        =
        \ell(4231)-\ell(1324).
    \]
    Hence this simultaneous embedding is an interval pattern embedding, and therefore
    \[
        [1324,4231]\cong [u,w].
    \]
\end{example}    

\subsection{Schubert varieties and singular loci}
\label{subsec:KL-singular-loci}

Let \(G\) be the general linear group \(GL_n(\mathbb C)\).
Let \(\mathcal B\subseteq G\) be the Borel subgroup of upper triangular matrices, and let \(T\subseteq \mathcal B\) be the diagonal torus.
For \(w\in S_n\), identifying \(w\) with the corresponding permutation matrix in \(G\), write
\[
    X_w^\circ:=\mathcal B w\mathcal B/\mathcal B,
    \qquad
    X_w:=\overline{X_w^\circ}
\]
for the Schubert cell and Schubert variety in the flag variety \(G/\mathcal B\).
For \(u\in S_n\), let
\[
    e_u:=u\mathcal B/\mathcal B
\]
denote the associated \(T\)-fixed point.
The closure relations are governed by the Bruhat order:
\[
    X_w=\bigsqcup_{u\leq w}X_u^\circ.
\]

For \(u\leq w\), Deodhar's
criterion~\cite[Section~2(iv) and Main Theorem]{Deo85} gives
\[
    P_{u,w}(q)=1
    \quad\Longleftrightarrow\quad
    X_w\text{ is smooth along }X_u^\circ.
\]
Consequently,
\[
    \operatorname{Sing}(X_w)
    =
    \bigsqcup_{
    \substack{u\leq w\\P_{u,w}(q)\neq1}}X_u^\circ.
\]
Let \(\operatorname{MaxSing}(w)\) be the set of maximal elements of
\(\{u\leq w:P_{u,w}(q)\neq1\}\)
in Bruhat order.
The irreducible components of \(\operatorname{Sing}(X_w)\) are precisely the Schubert varieties \(X_v\) with \(v\in\operatorname{MaxSing}(w)\)~\cite[p.~3916]{BW03}.
Thus, if \(X_w\) is singular, its singular locus is irreducible if and only if \(\operatorname{MaxSing}(w)=\{v\}\) for some \(v\), in which case \(\operatorname{Sing}(X_w)=X_v\).

\subsection{Height and reverse-complement symmetry}
\label{subsec:height-symmetry}

Let \(a<b<c<d\) be a $3412$-embedding in \(w\in S_n\).
We write
\[
    \alpha:=w(a),\qquad
    \beta:=w(b),\qquad
    \gamma:=w(c),\qquad
    \delta:=w(d),
\]
so that \(\gamma<\delta<\alpha<\beta\).
The \emph{height} and \emph{amplitude} of the embedding are \(\alpha-\delta\) and \(\beta-\gamma\), respectively.
We use the minimum height $h(w)$ defined in the introduction.

We will use reverse-complement symmetry to reduce the
configurations considered in the proofs. Let \(w_0\) be the longest element of \(S_n\).
\begin{definition}
\label{def:reverse-complement}
    The \emph{reverse-complement symmetry} is the map
    \[
    \theta\colon S_n\longrightarrow S_n,
    \qquad
    \theta(w):=w_0ww_0;
    \]
    see \cite[Section~5]{UW13}. Equivalently,
    \[
    \theta(w)(i)=n+1-w(n+1-i).
    \]
\end{definition}

Since \(\theta(s_i)=s_{n-i}\), the map \(\theta\) is an
automorphism of the Coxeter system and therefore preserves
length, Bruhat order, and Kazhdan--Lusztig polynomials.
Consequently,
\(\operatorname{MaxSing}(\theta(w))
=\theta(\operatorname{MaxSing}(w))\).

To describe this symmetry geometrically, represent \(w\in S_n\)
by its graph
\[
    \Gamma(w):=\{(i,w(i)):i\in[n]\},
\]
with positions increasing from left to right and values from bottom to top.
Then \(\theta\) rotates this graph by \(180^\circ\) about
\(\bigl((n+1)/2,(n+1)/2\bigr)\).
This preserves each of the patterns \(3412\) and \(4231\),
as well as the height and amplitude of \(3412\)-embeddings.
In particular, \(h(\theta(w))=h(w)\).

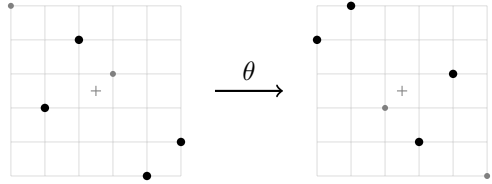
\begin{figure}[htbp]
    \centering
    \begin{tikzpicture}[scale=0.45,
        every node/.style={font=\small}]
        \begin{scope}
            \draw[gray!25, step=1] (1,1) grid (6,6);
            \draw[gray] (3.35,3.5) -- (3.65,3.5);
            \draw[gray] (3.5,3.35) -- (3.5,3.65);
            \foreach \x/\y in {1/6,4/4}
                \fill[gray] (\x,\y) circle (2.5pt);
            \foreach \x/\y in {2/3,3/5,5/1,6/2}
                \fill (\x,\y) circle (3.5pt);
        \end{scope}

        \draw[->, thick] (7,3.5) -- (9,3.5)
            node[midway,above] {\(\theta\)};

        \begin{scope}[xshift=9cm]
            \draw[gray!25, step=1] (1,1) grid (6,6);
            \draw[gray] (3.35,3.5) -- (3.65,3.5);
            \draw[gray] (3.5,3.35) -- (3.5,3.65);
            \foreach \x/\y in {3/3,6/1}
                \fill[gray] (\x,\y) circle (2.5pt);
            \foreach \x/\y in {1/5,2/6,4/2,5/4}
                \fill (\x,\y) circle (3.5pt);
        \end{scope}
    \end{tikzpicture}
    \caption{The highlighted \(3412\)-embedding and its image under reverse-complement.
    The crosses mark the center \((7/2,7/2)\) of the \(180^\circ\) rotation.}
    \label{fig:reverse-complement}
\end{figure}

\begin{example}
    Let \(w=635412\in S_6\), so that \(\theta(w)=563241\).
    The \(3412\)-embedding at positions \(2,3,5,6\) has height \(1\) and amplitude \(4\).
    Figure~\ref{fig:reverse-complement} shows its image under \(\theta\), which has the same height and amplitude.
\end{example}

\subsection{Singular locus components}
\label{subsec:singular-components}

To use the irreducibility hypothesis, we need criteria for
distinguishing singular locus components. We recall their
classification from \cite[Theorem~37]{BW03}, using the
interval-pattern formulation of \cite[Section~2.4]{Woo09}.
We then show that a $4231$-embedding guarantees a component
of type~I.

\begin{definition}
\label{def:types-of-irreducible-components}
    For a decreasing string, we use the notation
    \(j\cdots i=j(j-1)\cdots i\), with the string empty when \(j<i\).
    An irreducible component $X_u$ of $\operatorname{Sing}(X_w)$ 
    is said to be of one of the following types:
    \begin{enumerate}[label=(\roman*)]
        \item \emph{type~I} if, for some integers $y,z>0$, there exists an interval
        pattern embedding of 
        \begin{align*}
            \bigl[&
            (y+1)\cdots1\,(y+z+2)\cdots(y+2),\\
            &
            (y+z+2)(y+1)\cdots2\,
            (y+z+1)\cdots(y+2)1
            \bigr]
        \end{align*}
        into $[u,w]$.
        \item \emph{type~IIA} if, for some integers $y,z\geq0$, there exists an
        interval pattern embedding of
        \begin{align*}
            \bigl[&
            (y+1)\cdots1\,(y+3)(y+2)
            (y+z+4)\cdots(y+4),\\
            &
            (y+3)(y+1)\cdots2\,(y+z+4)1
            (y+z+3)\cdots(y+4)(y+2)
            \bigr]
        \end{align*}
        into $[u,w]$.
    
        \item \emph{type~IIB} if, for some integer $y\geq1$, there exists an interval
        pattern embedding of
        \[
            \bigl[
            1(y+3)\cdots2\,(y+4),\ 
            (y+3)(y+4)(y+2)\cdots312
            \bigr]
        \]
        into $[u,w]$.
    \end{enumerate}
    
    We refer to types~IIA and~IIB collectively as \emph{type~II}.
\end{definition}

Every irreducible component of \(\operatorname{Sing}(X_w)\) has exactly one of these types~\cite[Theorem~37]{BW03}.

\begin{lemma}
\label{lem:4231-gives-type-I}
    If $w$ contains $4231$, then the singular locus of $X_w$ has an irreducible component of type~I.
\end{lemma}

\begin{proof}
    Choose positions \(i<j<k<\ell\) forming a $4231$-embedding such that $\ell-i$ is minimal and,
    subject to this, such that $w(k)-w(j)$ is minimal.

    Suppose that this embedding does not index a type~I component.
    It is therefore a ``useless'' embedding in the terminology of~\cite[Section~13]{BW03}, so it must occur as one of the six dotted \(4231\)-embeddings in~\cite[Equation~(13.1)]{BW03}.
    The table below records the six possibilities and, in each case, another \(4231\)-embedding with positions \(i'<j'<k'<\ell'\):
    \[
        \begin{array}{c|c|c}
            \text{larger pattern} & \text{chosen embedding} & \text{new embedding}\\
            \hline
            52341  & 5241 & 5341\\
            52431  & 5241 & 5231\\
            53241  & 5241 & 5341\\
            53421  & 5341 & 5342\\
            54231  & 5231 & 4231\\
            635241 & 6341 & 5241
        \end{array}
    \]
    
    In the first three cases, \(\ell'-i'=\ell-i\), while \(w(k')-w(j')<w(k)-w(j)\); in the last three cases, \(\ell'-i'<\ell-i\).
    This contradicts the choice of the original embedding.
    Hence it indexes a type~I component.
\end{proof}

\subsection{Critical \texorpdfstring{$3412$}{3412}-embeddings}
\label{subsec:critical-3412-embeddings}

We now describe singular locus components arising from
$3412$-embeddings. The criteria below identify component types
and, in some configurations, give two distinct components.
Such configurations are excluded when the singular locus
is irreducible.

\begin{definition}
\label{def:critical-3412}
    Let \(a<b<c<d\) be a $3412$-embedding in $w$, with
    \[
        \alpha=w(a),\qquad \beta=w(b),\qquad
        \gamma=w(c),\qquad \delta=w(d),
    \]
    so that $\gamma<\delta<\alpha<\beta$.
    
    The \emph{critical regions} of the embedding are the four open rectangles
    \[
        (a,b)\times(\delta,\alpha),\qquad
        (b,c)\times(\alpha,\beta),\qquad
        (c,d)\times(\delta,\alpha),\qquad
        (b,c)\times(\gamma,\delta).
    \]
    The embedding is \emph{critical} if none of these regions contains a
    point of $\Gamma(w)$.
    
    For a critical embedding, define
    \[
        B:=\{p:b<p<c,\ \delta<w(p)<\alpha\},
    \]
    \[
        A_1:=\{p:a<p<b,\ \gamma<w(p)<\delta\},
    \]
    \[
        A_2:=\{p:c<p<d,\ \alpha<w(p)<\beta\},
    \]
    and set \(A:=A_1\cup A_2\).
    The embedding is \emph{reduced} if the values $w(p)$ for $p\in B$
    decrease from left to right.
\end{definition}

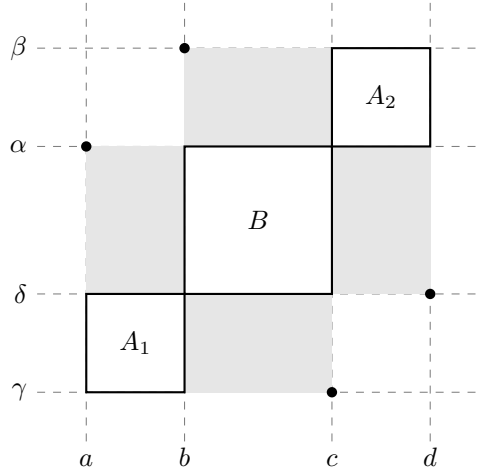
\begin{figure}[htbp]
    \centering
    \begin{tikzpicture}[scale=0.65, every node/.style={font=\small}]
        \def\xa{1}
        \def\xb{3}
        \def\xc{6}
        \def\xd{8}
    
        \def\ygamma{1}
        \def\ydelta{3}
        \def\yalpha{6}
        \def\ybeta{8}
    
        \draw[dashed, gray] (\xa,0) -- (\xa,9);
        \draw[dashed, gray] (\xb,0) -- (\xb,9);
        \draw[dashed, gray] (\xc,0) -- (\xc,9);
        \draw[dashed, gray] (\xd,0) -- (\xd,9);
    
        \draw[dashed, gray] (0,\ygamma) -- (9,\ygamma);
        \draw[dashed, gray] (0,\ydelta) -- (9,\ydelta);
        \draw[dashed, gray] (0,\yalpha) -- (9,\yalpha);
        \draw[dashed, gray] (0,\ybeta) -- (9,\ybeta);
    
        \fill[gray!20] (\xa,\ydelta) rectangle (\xb,\yalpha);
        \fill[gray!20] (\xb,\yalpha) rectangle (\xc,\ybeta);
        \fill[gray!20] (\xc,\ydelta) rectangle (\xd,\yalpha);
        \fill[gray!20] (\xb,\ygamma) rectangle (\xc,\ydelta);
    
        \draw[thick] (\xa,\ygamma) rectangle (\xb,\ydelta);
        \draw[thick] (\xb,\ydelta) rectangle (\xc,\yalpha);
        \draw[thick] (\xc,\yalpha) rectangle (\xd,\ybeta);
    
        \node at (2,2) {$A_1$};
        \node at (4.5,4.5) {$B$};
        \node at (7,7) {$A_2$};
    
        \fill (\xa,\yalpha) circle (3pt);
        \fill (\xb,\ybeta) circle (3pt);
        \fill (\xc,\ygamma) circle (3pt);
        \fill (\xd,\ydelta) circle (3pt);
    
        \node[anchor=base] at (\xa,-0.5) {$a$};
        \node[anchor=base] at (\xb,-0.5) {$b$};
        \node[anchor=base] at (\xc,-0.5) {$c$};
        \node[anchor=base] at (\xd,-0.5) {$d$};
    
        \node[left] at (0,\ygamma) {$\gamma$};
        \node[left] at (0,\ydelta) {$\delta$};
        \node[left] at (0,\yalpha) {$\alpha$};
        \node[left] at (0,\ybeta) {$\beta$};
    \end{tikzpicture}
    \caption{The regions associated with a $3412$-embedding $a<b<c<d$.
    The four shaded rectangles are the critical regions; the unshaded ones correspond to $A_1$, $B$, and $A_2$.}
    \label{fig:3412-regions}
\end{figure}

By abuse of terminology, we refer to $A_1,A_2$, and $B$ as regions; these sets record the positions of the points of $\Gamma(w)$ in the corresponding rectangles.
Figure~\ref{fig:3412-regions} illustrates these regions associated with a $3412$-embedding.

The next lemma connects these regions with the minimum height
$h(w)$ appearing in the main theorem.

\begin{lemma}
\label{lem:minimal-3412-critical-reduced}
    Let $a<b<c<d$ be a $3412$-embedding in $w$ of minimum height and, among those, minimum amplitude.
    Then the embedding is critical and reduced.
    Moreover, writing $\alpha=w(a)$ and $\delta=w(d)$, its associated region $B$ satisfies
    \[
        \{w(p):p\in B\}=\{\delta+1,\ldots,\alpha-1\}.
    \]
    In particular, $|B|=h(w)-1$.
\end{lemma}

\begin{proof}
    Any point in a critical region would give a $3412$-embedding of
    smaller height, or of the same height and smaller amplitude.
    Hence the embedding is critical.
    If $p<q$ lie in $B$ with $w(p)<w(q)$, then $a<b<p<q$ is a $3412$-embedding of smaller height.
    Thus the embedding is reduced.

    Let $t$ be an integer with $\delta<t<\alpha$, and set $p:=w^{-1}(t)$.
    If $p<b$, then $p<b<c<d$ is a $3412$-embedding of height $t-\delta<\alpha-\delta$.
    If $p>c$, then $a<b<c<p$ is a $3412$-embedding of height $\alpha-t<\alpha-\delta$.
    Both contradict the minimality of the height.
    Therefore $b<p<c$, so $p\in B$.
    This proves the claimed equality and hence $|B|=\alpha-\delta-1=h(w)-1$.
\end{proof}

For a reduced critical embedding, the following criteria
relate the regions $A_1,A_2,B$ to singular locus components.

\begin{proposition}
[{\cite[Section~2.4]{Woo09}}]
\label{prop:woo-component-dictionary}
    Let $a<b<c<d$ be a reduced critical $3412$-embedding in $w$, with
    associated regions $A,A_1,A_2$, and $B$ as above.
    \begin{enumerate}[label=(\roman*), ref=(\roman*)]
        \item\label{item:woo-IIA}
        If $B=\varnothing$, there exists $v\in\operatorname{MaxSing}(w)$ such that $X_v$ is of type~IIA.
    
        \item\label{item:woo-IIB}
        If $A=\varnothing$ and $B\neq\varnothing$, there exists $v\in\operatorname{MaxSing}(w)$ such that $X_v$ is of type~IIB.
    
        \item\label{item:woo-I-left}\label{item:woo-I-right}
        If $A_i\neq\varnothing$ for some \(i\in\{1,2\}\) and $B\neq\varnothing$, there exists $v\in\operatorname{MaxSing}(w)$ such that $X_v$ is of type~I.
    
        \item\label{item:woo-two-I}
        If $A_1,A_2$, and $B$ are all nonempty, there exist distinct $v_1,v_2\in\operatorname{MaxSing}(w)$ such that both $X_{v_1}$ and $X_{v_2}$ are of type~I.
    \end{enumerate}
\end{proposition}

\subsection{Restrictions and Kazhdan--Lusztig inequalities}
\label{subsec:pattern-maps-KL-inequalities}

We first record the following consequence of
\cite[Theorem~1]{BBr03}.

\begin{lemma}
\label{lem:pattern-lower-bound}
    If a permutation $w$ contains \(y\), then \(P_{\mathrm{id},y}(1)\leq P_{\mathrm{id},w}(1)\).
\end{lemma}

We now turn to coefficientwise inequalities for Kazhdan--Lusztig polynomials.
Write $[q^j]f(q)$ for the coefficient of $q^j$ in $f(q)$.
For \(f,g\in\mathbb Z[q]\), write \(f\leq_{\mathrm{coeff}}g\) if every coefficient of \(g-f\) is nonnegative.

\begin{theorem}[{\cite[Corollary~4]{Ir88}}]
\label{thm:KL-monotonicity}
    For \(u\le v\le w\) in \(S_n\), we have 
    \[
        P_{v,w}(q)\leq_{\mathrm{coeff}}P_{u,w}(q).
    \]
\end{theorem}   

\begin{definition}
\label{def:interval-restrictions}
    Let $w\in S_n$ and let $J=[a,b]\subseteq[n]$ be a nonempty interval.
    Let $i_1<\cdots<i_k$ be the positions whose values lie in $J$, where $k=b-a+1$.
    The \emph{value restriction} and \emph{position restriction} of $w$ to $J$ are, respectively,
    \[
        \begin{aligned}
            w|_J&:=\operatorname{fl}\bigl(w(i_1)\cdots w(i_k)\bigr),\\
            w|^J&:=\operatorname{fl}\bigl(w(a)w(a+1)\cdots w(b)\bigr).
        \end{aligned}
    \]
\end{definition}

\begin{proposition}
\label{prop:continuous-restriction}
    For every $w\in S_n$ and every nonempty interval $J\subseteq[n]$,
    \[
        P_{\mathrm{id},w|_J}(q)
        \leq_{\mathrm{coeff}}
        P_{\mathrm{id},w}(q),
        \qquad
        P_{\mathrm{id},w|^J}(q)
        \leq_{\mathrm{coeff}}
        P_{\mathrm{id},w}(q).
    \]
\end{proposition}

\begin{proof}
    Write $J=[a,b]$, set \(W':=\langle s_a,\ldots,s_{b-1}\rangle\), and let $x$ be the minimal-length element of the coset $W'w$.
    Under the identification \(W'\simeq S_{b-a+1}\), 
    restriction to the values in $J$
    sends $x$ to $\mathrm{id}$ and $w$ to $w|_J$.
    Hence \cite[Theorem~6]{BBr03} gives
    \[
        P_{x,w}(q)=P_{\mathrm{id},w|_J}(q).
    \]
    Since $\mathrm{id}\leq x\leq w$, the first inequality follows from Theorem~\ref{thm:KL-monotonicity}.
    
    The second follows by applying the first to $w^{-1}$ and using
    \[
        (w|^J)^{-1}=(w^{-1})|_J,
        \qquad
        P_{u,w}(q)=P_{u^{-1},w^{-1}}(q).\qedhere
    \]
\end{proof}

Finally, we record the following consequence of
\cite[Theorems~2, 5, and~6]{BBr03}.

\begin{theorem}
\label{thm:Billey-Braden-product}
    Let $W'\subseteq S_n$ be a standard parabolic subgroup and let $w\in S_n$.
    Then $[\mathrm{id},w]\cap W'$ has a greatest element
    $y$ in Bruhat order, and
    \[
        P_{\mathrm{id},y}(q)
        \leq_{\mathrm{coeff}}
        P_{\mathrm{id},w}(q).
    \]
\end{theorem}

%% file: 3_proof_R1.tex
\section{The 3412-only case}
\label{sec:3412-only}

In this section, we prove~\ref{case:R1}.
Let $w\in S_n$ satisfy the hypotheses of Theorem~\ref{thm:main}.
Assume that $w$ contains $3412$ and avoids $4231$, and set $h=h(w)$.

We follow the construction and notation of \cite[Section~4.1]{Woo09}.
Choose a $3412$-embedding $a<b<c<d$ of minimum height and, among those, minimum amplitude.
Write
\[
    \alpha:=w(a),\qquad
    \beta:=w(b),\qquad
    \gamma:=w(c),\qquad
    \delta:=w(d),
\]
so that $\gamma<\delta<\alpha<\beta$ and $\alpha-\delta=h$.
By Lemma~\ref{lem:minimal-3412-critical-reduced}, this embedding
is critical and reduced.
Define
\begin{gather*}
    \alpha'
    :=
    \max\{j\in[n]:j\geq\alpha,\,
    w^{-1}(j)<w^{-1}(j-1)<\cdots<w^{-1}(\alpha)\},\\
    \delta'
    :=
    \min\{j\in[n]:j\leq\delta,\,
    w^{-1}(\delta)<w^{-1}(\delta-1)<\cdots<w^{-1}(j)\},
\end{gather*}
and set
\[
    \kappa:=\delta'+\alpha'-\alpha.
\]
Then \(\gamma<\delta'\leq\delta<\alpha\leq\alpha'<\beta\).
Let
\[
    I:=\{s_{\delta'},\ldots,s_{\alpha'-1}\},
    \qquad
    J:=I\setminus\{s_\kappa\},
    \qquad
    v:=w_Jw_Iw,
\]
where $w_I$ and $w_J$ are the longest elements of the standard parabolic subgroups of $S_n$ generated by $I$ and $J$, respectively.
Consider the variety
\[
    Z:=P_I\times^{P_J}X_v
\]
and the proper birational morphism
\[
    \pi\colon Z\longrightarrow X_w,
    \qquad [p,F_\bullet]\longmapsto p\cdot F_\bullet,
\]
where $P_I,P_J$ are the corresponding standard
parabolic subgroups of $\mathrm{GL}_n$.
Since $w$ avoids $4231$, it also avoids $463152$.
Together with the irreducibility of the singular locus, \cite[Lemma~4.1]{Woo09} therefore shows that $Z$ is smooth.
Thus $\pi$ is a resolution of singularities, known as \emph{Cortez's resolution}.

\begin{figure}[htbp]
    \centering
    \begin{tikzpicture}[x=0.60cm,y=0.60cm,font=\small]

        \def\xmin{0.0}
        \def\xmax{13.2}
        \def\ymin{0.1}
        \def\ymax{10.7}

        \def\xa{3.0}
        \def\xb{5.4}
        \def\xc{8.0}
        \def\xd{10.4}

        \def\xap{1.0}
        \def\xdp{12.3}

        \def\yg{1.6}
        \def\ydp{3.0}
        \def\yd{4.3}
        \def\ya{7.0}
        \def\yap{8.2}
        \def\yb{9.4}

        \tikzset{
            LRregion/.style={fill=gray!20},
            emptyregion/.style={
                fill=blue!8,
                pattern=north east lines,
                pattern color=blue!60
            }
        }

        \fill[LRregion] (\xa,\yb) rectangle (\xb,\ymax); 
        \fill[LRregion] (\xc,\ymin) rectangle (\xd,\yg); 


        \fill[emptyregion] (\xmin,\yap) rectangle (\xa,\ymax);

        \fill[emptyregion] (\xa,\yap) rectangle (\xb,\yb);

        \fill[emptyregion] (\xb,\yap) rectangle (\xc,\ymax);

        \fill[emptyregion] (\xc,\yg) rectangle (\xd,\ydp);

        \fill[emptyregion] (\xd,\ymin) rectangle (\xmax,\ydp);


        \draw[dashed, gray] (\xa,\ya) -- (\xa,\ymax);
        \draw[dashed, gray] (\xb,\ymin) -- (\xb,\ymax);
        \draw[dashed, gray] (\xc,\ymin) -- (\xc,\ymax);
        \draw[dashed, gray] (\xd,\ymin) -- (\xd,\yd);

        \draw[dashed, gray] (\xmin,\yb) -- (\xc,\yb);
        \draw[dashed, gray] (\xmin,\ya) -- (\xmax,\ya);
        \draw[dashed, gray] (\xmin,\yd) -- (\xmax,\yd);
        \draw[dashed, gray] (\xb,\yg) -- (\xmax,\yg);

        \draw[dashed, gray] (\xmin,\yap) -- (\xc,\yap);
        \draw[dashed, gray] (\xc,\ydp) -- (\xmax,\ydp);


        \draw[densely dotted, thick] (\xap,\yap) -- (\xa,\ya);
        \fill (1.65,7.80) circle (1.4pt);
        \fill (2.30,7.38) circle (1.4pt);

        \draw[densely dotted, thick] (\xd,\yd) -- (\xdp,\ydp);
        \fill (11.05,3.86) circle (1.4pt);
        \fill (11.70,3.43) circle (1.4pt);

        \fill (\xa,\ya) circle (2.6pt);
        \fill (\xb,\yb) circle (2.6pt);
        \fill (\xc,\yg) circle (2.6pt);
        \fill (\xd,\yd) circle (2.6pt);
        \fill (\xap,\yap) circle (2.2pt);
        \fill (\xdp,\ydp) circle (2.2pt);

        \node[below left]  at (\xa,\ya) {$a$};
        \node[below right] at (\xb,\yb) {$b$};
        \node[below left]  at (\xc,\yg) {$c$};
        \node[below left]  at (\xd,\yd) {$d$};

        \node at (4.20,10.05) {$R$};
        \node at (9.20,0.82) {$L$};

        \node[left]  at (\xmin,\yb)  {$\beta$};
        \node[left]  at (\xmin,\yap) {$\alpha'$};
        \node[left]  at (\xmin,\ya)  {$\alpha$};

        \node[right] at (\xmax,\yd)  {$\delta$};
        \node[right] at (\xmax,\ydp) {$\delta'$};
        \node[right] at (\xmax,\yg)  {$\gamma$};

    \end{tikzpicture}
    \caption{The chosen $3412$-embedding.
    The gray shaded regions correspond to $R$ and $L$.
    The dotted segments indicate the consecutive strings from
    $\alpha'$ to $\alpha$ and from $\delta$ to $\delta'$.
    The blue hatched regions are empty by
    \cite[Lemma~5.1, Figure~4]{Woo09}.}
    \label{fig:R1-identity-configuration}
\end{figure}
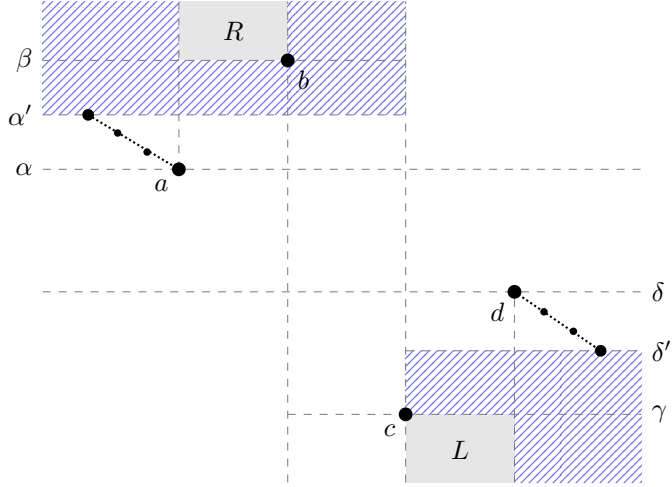

To determine the identity fiber explicitly, let
\[
    r_+:=\alpha'-\alpha,
    \qquad
    r_-:=\delta-\delta',
\]
and define
\[
    L:=\{p:c<p<d,\ w(p)<\gamma\},
    \qquad
    R:=\{p:a<p<b,\ w(p)>\beta\}.
\]
By abuse of terminology, we also refer to $L$ and $R$ as regions; these sets record the positions of the points of $\Gamma(w)$ in the corresponding rectangles.
Figure~\ref{fig:R1-identity-configuration} illustrates the regions $L,R$ and the sequences defining $\alpha'$ and $\delta'$.

\begin{proposition}
\label{prop:R1-identity-fiber}
    The fiber over $e_{\mathrm{id}}$ is
    \begin{equation}
    \label{eq:R1-id-fiber}
        \pi^{-1}(e_{\mathrm{id}})
        \cong
        \operatorname{Gr}_x(\mathbb C^{x+y}),
    \end{equation}
    where $x=1+\min\{r_+,|L|\}$ and
    $y=h+\min\{r_-,|R|\}$.
\end{proposition}

\begin{proof}
    Set
    \[
        M:=\max\bigl(\{p<c:w(p)<\delta'\}\cup\{a\}\bigr),
        \qquad
        N:=\max\{p:w(p)<\delta'\}.
    \]
    Applying \cite[Lemma~4.2]{Woo09} to $e_{\mathrm{id}}$ gives \eqref{eq:R1-id-fiber}, with
    \[
        x=\kappa-\max\{\delta'-1,M\},
        \qquad
        y=\min\{\alpha',N\}-\kappa.
    \]
    Since $\kappa-(\delta'-1)=1+r_+$ and $\alpha'-\kappa=h+r_-$, it remains to prove
    \begin{gather}
        \kappa-M=1+|L|,
        \label{eq:R1-kappa-M}\\
        N-\kappa=h+|R|.
        \label{eq:R1-N-kappa}
    \end{gather}

    For $\sigma\in S_n$, define its \emph{rank function} by
    \[
        r_\sigma(p,t):=\#\{i\leq p:\sigma(i)\leq t\}.
    \]
    As shown in the proof of \cite[Lemma~4.2]{Woo09},
    \[
        r_v(M,\kappa)=M,
        \qquad
        r_v(N,\kappa)=\kappa.
    \]
    Since $\kappa=\delta'+\alpha'-\alpha$, the action of $w_Jw_I$ on the values of $w$ is
    \[
        v(p)=
        \begin{cases}
            w(p),
            & w(p)<\delta' \text{ or } w(p)>\alpha',\\
            w(p)+\alpha'-\alpha+1,
            & \delta'\leq w(p)<\alpha,\\
            w(p)-(\alpha-\delta'),
            & \alpha\leq w(p)\leq\alpha'.
        \end{cases}
    \]
    Hence
    \begin{align}
        v(p)\leq\kappa
        \quad&\Longleftrightarrow\quad
        w(p)<\delta'
        \text{ or }
        \alpha\leq w(p)\leq\alpha',
        \label{eq:R1-v-le-kappa}\\
        v(p)>\kappa
        \quad&\Longleftrightarrow\quad
        \delta'\leq w(p)<\alpha
        \text{ or }
        w(p)>\alpha'.
        \label{eq:R1-v-gt-kappa}
    \end{align}

    We claim that
    \begin{align*}
        \{p>M:w(p)<\delta'\}&=\{c\}\sqcup L,\\
        \{p<N:w(p)>\alpha'\}&=\{b\}\sqcup R.
    \end{align*}
    For the first identity, the definition of $M$ excludes positions $M<p<c$ with $w(p)<\delta'$, while $w(c)=\gamma<\delta'$.
    The empty regions illustrated in Figure~\ref{fig:R1-identity-configuration} show that the remaining positions are precisely those in $L$.
    Since $M<c$ and every element of $L$ lies between $c$ and $d$, this also gives $c\leq N<d$.

    For the second identity, the same figure shows that the positions $p\leq c$ with $w(p)>\alpha'$ are precisely $b$ and those in $R$.
    There are no additional such positions between $c$ and $N$ by \cite[Lemma~5.2(iv)]{Woo09}.

    We first compute $\kappa-M$.
    Since $v$ has exactly $\kappa$ values at most $\kappa$ and $r_v(M,\kappa)=M$,
    \[
        \kappa-M=\#\{p>M:v(p)\leq\kappa\}.
    \]
    By the definition of $\alpha'$, $w^{-1}(\alpha')<\cdots<w^{-1}(\alpha)=a\leq M$, so no value in $[\alpha,\alpha']$ occurs after $M$.
    Using \eqref{eq:R1-v-le-kappa} and the first set identity gives
    \[
        \kappa-M
        =\#\{p>M:w(p)<\delta'\}
        =1+|L|,
    \]
    proving \eqref{eq:R1-kappa-M}.

    We next compute $N-\kappa$.
    Since $w(N)<\delta'$, \eqref{eq:R1-v-le-kappa} gives $v(N)\leq\kappa$.
    Thus $r_v(N,\kappa)=\kappa$ and
    \eqref{eq:R1-v-gt-kappa} yield
    \begin{align*}
        N-\kappa
        &=\#\{p<N:v(p)>\kappa\}\\
        &=\#\{p<N:\delta'\leq w(p)<\alpha\}
          +\#\{p<N:w(p)>\alpha'\}.
    \end{align*}
    Lemma~\ref{lem:minimal-3412-critical-reduced} gives $b<w^{-1}(t)<c\leq N$ for $\delta<t<\alpha$, whereas the definition of $\delta'$ gives $w^{-1}(t)\geq d>N$ for $\delta'\leq t\leq\delta$.
    Hence
    \[
        \#\{p<N:\delta'\leq w(p)<\alpha\}=h-1.
    \]
    The second set identity gives $\#\{p<N:w(p)>\alpha'\}=1+|R|$.
    Consequently,
    \[
        N-\kappa=(h-1)+(1+|R|)=h+|R|,
    \]
    proving \eqref{eq:R1-N-kappa} and completing the proof.
\end{proof}

We next record how $x$ and $y$ transform under reverse-complement.

\begin{corollary}
\label{cor:R1-reverse-complement-fiber}
    Let $\widehat x,\widehat y$ be the integers in
    Proposition~\ref{prop:R1-identity-fiber} for the image of the chosen
    embedding in $\theta(w)$.
    Then
    \[
        \widehat x=1+y-h,
        \qquad
        \widehat y=h+x-1.
    \]
\end{corollary}

\begin{proof}
    For this embedding, $\widehat\alpha'=n+1-\delta'$ and $\widehat\delta'=n+1-\alpha'$.
    Hence
    $\widehat r_+=r_-$ and $\widehat r_-=r_+$.
    Moreover, the $180^\circ$ rotation of the permutation graph gives
    $|\widehat L|=|R|$ and $|\widehat R|=|L|$,
    and $\widehat h=h$ by reverse-complement symmetry.
    Therefore, by Proposition~\ref{prop:R1-identity-fiber},
    \begin{gather*}
        \widehat x
        =
        1+\min\{r_-,|R|\}
        =
        1+y-h,\\
        \widehat y
        =
        h+\min\{r_+,|L|\}
        =
        h+x-1.\qedhere
    \end{gather*}
\end{proof}

We now combine the formulas above with the hypothesis
$P_{\mathrm{id},w}(1)=3$ to restrict the possible values of $x$ and $y$.

\begin{lemma}
\label{lem:R1-fiber-bounds}
We have
\[
    x\leq 2,
    \qquad
    y\leq h+1,
    \qquad
    (x,y)\neq(2,h+1).
\]
\end{lemma}

\begin{proof}
    The Kazhdan--Lusztig recurrence gives
    \[
        P_{\mathrm{id},65872143}(1)=6,\qquad
        P_{\mathrm{id},57861432}(1)=5,\qquad
        P_{\mathrm{id},58761432}(1)=4.
    \]
    
    Suppose first that $y\geq h+2$.
    By Proposition~\ref{prop:R1-identity-fiber}, we have
    $r_-\geq2$ and $|R|\geq2$.
    Choose $f_1<f_2$ in $R$, and set $d_1:=w^{-1}(\delta-1)$ and $d_2:=w^{-1}(\delta-2)$.
    Since $r_-\geq2$, the definition of $\delta'$ gives
    $d<d_1<d_2$.
    Thus
    \[
        a<f_1<f_2<b<c<d<d_1<d_2.
    \]
    The entries in these positions form the pattern $57861432$ or
    $58761432$, according to the relative order of $w(f_1)$ and $w(f_2)$.
    Lemma~\ref{lem:pattern-lower-bound} then gives $P_{\mathrm{id},w}(1)\geq4$, a contradiction.
    Hence $y\leq h+1$.
    
    By reverse-complement symmetry, the same bound applied to \(\theta(w)\) gives \(\widehat y\leq h+1\).
    By Corollary~\ref{cor:R1-reverse-complement-fiber}, \(\widehat y=h+x-1\), and therefore \(x\leq2\).
    
    It remains to exclude $(x,y)=(2,h+1)$.
    By Proposition~\ref{prop:R1-identity-fiber}, this would imply $r_+,r_-\geq1$ and $L,R\neq\varnothing$.
    Choose $f\in R$ and $g\in L$, and set $a_1:=w^{-1}(\alpha+1)$ and $d_1:=w^{-1}(\delta-1)$.
    The definitions of $\alpha'$ and $\delta'$ give $a_1<a$ and $d<d_1$.
    Hence
    \[
        a_1<a<f<b<c<g<d<d_1.
    \]
    The entries in these positions form the pattern $65872143$.
    Lemma~\ref{lem:pattern-lower-bound} would then give
    $P_{\mathrm{id},w}(1)\geq6$, again a contradiction.
\end{proof}

Since Proposition~\ref{prop:R1-identity-fiber} also gives $x\geq 1$ and
$y\geq h$, Lemma~\ref{lem:R1-fiber-bounds} leaves exactly three
possibilities:
\begin{equation}
\label{eq:R1-fiber-possibilities}
    (x,y)\in\{(1,h),(1,h+1),(2,h)\}.
\end{equation}
We now combine these possibilities with the decomposition theorem and Theorem~\ref{thm:Gaetz-Gao-min-degree} to determine $P_{\mathrm{id},w}(q)$.
For $z\leq w$, set
\[
    H_{z,\pi}(q)
    :=
    \sum_{i\geq0}
    q^i\dim H^{2i}\bigl(\pi^{-1}(e_z);\mathbb C\bigr).
\]
Applying the decomposition theorem as in \cite[Section~4.3]{Woo09} gives the following consequence.

\begin{proposition}
\label{prop:R1-decomposition}
    For the resolution $\pi\colon Z\to X_w$ constructed above, there exist polynomials $D_u(q)\in\mathbb Z_{\geq0}[q]$ for $u<w$ such
    that, for every $z\leq w$,
    \[
        H_{z,\pi}(q)
        =
        P_{z,w}(q)
        +
        \sum_{z\leq u<w}
        D_u(q)P_{z,u}(q).
    \]
\end{proposition}

\begin{proof}[Proof of~\ref{case:R1}]
    By \eqref{eq:R1-fiber-possibilities} and
    Corollary~\ref{cor:R1-reverse-complement-fiber}, after replacing $w$ by its reverse-complement if necessary, we may assume $x=1$ and $y\in\{h,h+1\}$.
    By Proposition~\ref{prop:R1-identity-fiber},
    $\pi^{-1}(e_{\mathrm{id}})\cong\mathbb P^y$.
    Therefore Proposition~\ref{prop:R1-decomposition} gives
    \[
        P_{\mathrm{id},w}(q)
        \leq_{\mathrm{coeff}}
        H_{\mathrm{id},\pi}(q)
        =
        1+q+\cdots+q^y
        \leq_{\mathrm{coeff}}
        1+q+\cdots+q^{h+1}.
    \]
    
    Since $w$ avoids $4231$,
    Theorem~\ref{thm:Gaetz-Gao-min-degree} gives
    $[q^j]P_{\mathrm{id},w}(q)=0$ for $1\leq j<h$.
    Hence
    \[
        P_{\mathrm{id},w}(q)
        =
        1+\varepsilon_1 q^h+\varepsilon_2 q^{h+1},
        \qquad
        \varepsilon_1,\varepsilon_2\in\{0,1\}.
    \]
    Finally, $P_{\mathrm{id},w}(1)=3$ forces
    $\varepsilon_1=\varepsilon_2=1$.
\end{proof}

\begin{remark}
    The possibility $(x,y)=(1,h)$ is in fact excluded by
    $P_{\mathrm{id},w}(1)=3$.
    Indeed, the same argument as above gives \(P_{\mathrm{id},w}(q)=1+q^h\).
    This is the identity-fiber calculation used in
    \cite[Section~4]{Woo09}.
\end{remark}

%% file: 4_proof_R2.tex
\section{The 4231-only case}
\label{sec:4231-only}

In this section, we prove~\ref{case:R2}.
Let $w\in S_n$ satisfy the hypotheses of Theorem~\ref{thm:main}.
Assume that $w$ avoids $3412$ and contains $4231$.

A permutation is called \emph{covexillary} if it avoids $3412$.
Lascoux gave a formula for the Kazhdan--Lusztig polynomials of
covexillary permutations~\cite{Las95}.
We use its rooted-tree formulation (Theorem~\ref{thm:R2-Lascoux-tree}).

For a covexillary permutation $x$, the rooted-tree formulation
associates a rooted tree $T_x$ to $x$. In the case of
$P_{\mathrm{id},x}(q)$, each leaf $\lambda$ of $T_x$ is assigned a
nonnegative integer \emph{capacity} $c(\lambda)$.
Write $E(T_x)$ for the edge set of $T_x$.

\begin{definition}
\label{def:R2-admissible-labeling}
An \emph{edge labeling} of $T_x$ is a map
\(L:E(T_x)\longrightarrow \mathbb Z_{\geq 0}\).
Set
\[
    s(L):=\sum_{e\in E(T_x)}L(e).
\]
We call $L$ \emph{admissible} if, for every leaf $\lambda$, whenever
$e_1,\ldots,e_r$ are the edges on the path from the root to $\lambda$
in this order, one has
\[
    L(e_1)\leq\cdots\leq L(e_r)\leq c(\lambda).
\]
\end{definition}

\begin{theorem}[{\cite[Section~3]{WooFPSAC09}}]
\label{thm:R2-Lascoux-tree}
For every covexillary permutation $x$,
\[
    P_{\mathrm{id},x}(q)
    =
    \sum_L q^{s(L)},
\]
where the sum ranges over all admissible edge labelings of $T_x$.
\end{theorem}

We now use this formula to prove~\ref{case:R2}.

\begin{lemma}
\label{lem:R2-terminal-chain}
    There exist integers $c,d\geq1$ such that
    \[
        P_{\mathrm{id},w}(q)
        =
        \sum_{\substack{a_1,\ldots,a_d\in\mathbb Z\\
                        0\leq a_1\leq\cdots\leq a_d\leq c}}
        q^{a_1+\cdots+a_d}.
    \]
\end{lemma}

\begin{proof}
    By \cite[Section~3]{WooFPSAC09}, the condition
    $|\operatorname{MaxSing}(w)|=1$ implies that $T_w$ has a unique leaf
    $\lambda$ of positive capacity $c:=c(\lambda)>0$.
    
    Let $\mu$ be the last branching vertex on the path from the root to
    $\lambda$, or the root if no such vertex exists.
    Let $d$ be the number of edges from $\mu$ to $\lambda$.
    Since every edge outside this chain lies on a path to a leaf other than
    $\lambda$, whose capacity is zero, admissibility forces all such edges to be labeled zero.
    If $d=0$, Theorem~\ref{thm:R2-Lascoux-tree} would give
    $P_{\mathrm{id},w}(q)=1$, contradicting the singularity of $w$.
    Thus $d\geq1$.
    
    Let $a_1,\ldots,a_d$ be the edge labels on the chain from $\mu$ to $\lambda$. Admissibility is then equivalent to
    \(0\leq a_1\leq\cdots\leq a_d\leq c\),
    and every such tuple extends uniquely to an admissible edge labeling
    of $T_w$. Therefore, Theorem~\ref{thm:R2-Lascoux-tree} gives
    \[
        P_{\mathrm{id},w}(q)
        =
        \sum_{\substack{a_1,\ldots,a_d\in\mathbb Z\\
                        0\leq a_1\leq\cdots\leq a_d\leq c}}
        q^{a_1+\cdots+a_d}.\qedhere
    \]
\end{proof}

\begin{proof}[Proof of~\ref{case:R2}]
    Evaluating the formula in Lemma~\ref{lem:R2-terminal-chain} at $q=1$ gives
    \[
        3=P_{\mathrm{id},w}(1)=\binom{c+d}{d},
        \qquad c,d\geq1.
    \]
    Hence \((c,d)=(2,1)\) or \((1,2)\).
    In either case, Lemma~\ref{lem:R2-terminal-chain} gives
    \[
        P_{\mathrm{id},w}(q)=1+q+q^2.\qedhere
    \]
\end{proof}

%% file: 5_proof_R3.tex
\section{The mixed case}
\label{sec:mixed}

In this section, we prove~\ref{case:R3}.
Let $w\in S_n$ satisfy the hypotheses of Theorem~\ref{thm:main}.
Assume that $w$ contains both $3412$ and $4231$, and set $h=h(w)$.

\begin{lemma}
\label{lem:mixed-height-at-least-two}
    We have $h\geq2$.
\end{lemma}

\begin{proof}
    Suppose that $h=1$, and choose a $3412$-embedding of minimum
    height and, among those, minimum amplitude.
    By Lemma~\ref{lem:minimal-3412-critical-reduced}, this embedding
    is critical and reduced, with $B=\varnothing$.
    Proposition~\ref{prop:woo-component-dictionary} therefore gives
    an irreducible component of the singular locus of type~IIA.
    Since $w$ contains $4231$, Lemma~\ref{lem:4231-gives-type-I}
    also gives a component of type~I.
    These are distinct components, contradicting $|\operatorname{MaxSing}(w)|=1$.
\end{proof}

We now prove that $[q^h]P_{\mathrm{id},w}(q)>0$.
Choose a $3412$-embedding $a<b<c<d$ of minimum height
and, among those, minimum amplitude.
Write
\[
    \alpha:=w(a),\qquad
    \beta:=w(b),\qquad
    \gamma:=w(c),\qquad
    \delta:=w(d),
\]
so that $\gamma<\delta<\alpha<\beta$ and $\alpha-\delta=h$.
Let \(i_1<\cdots<i_N\) be the positions $p\in[a,d]$ such that
\(\gamma\leq w(p)\leq\beta\), and set
\[
    \rho:=\operatorname{fl}\bigl(w(i_1)\cdots w(i_N)\bigr)\in S_N.
\]

\begin{lemma}
\label{lem:R3-continuous-reduction}
We have
\[
    P_{\mathrm{id},\rho}(q)
    \leq_{\mathrm{coeff}}
    P_{\mathrm{id},w}(q)
\]
and $N\geq h+3$. Moreover, after replacing $\rho$ by its
reverse-complement if necessary,
\[
    \rho(1)\cdots\rho(h+2)
    =
    (h+2)\,N\,(h+1)\cdots3\,1.
\]
\end{lemma}

\begin{proof}
    By construction, we obtain $\rho$ by restricting $w$ to the values
    in $[\gamma,\beta]$ and then retaining the positions corresponding
    to $[a,d]$. These positions form an interval in the restricted
    permutation. Hence
    Proposition~\ref{prop:continuous-restriction} gives
    \[
        P_{\mathrm{id},\rho}(q)
        \leq_{\mathrm{coeff}}
        P_{\mathrm{id},w}(q).
    \]
    
    Let $A_1,A_2,B$ be the associated regions defined in
    Definition~\ref{def:critical-3412}.
    By Lemma~\ref{lem:minimal-3412-critical-reduced}, the chosen
    $3412$-embedding is critical and reduced, and the values occurring
    in $B$ are $\alpha-1,\ldots,\delta+1$, in decreasing order
    from left to right. In particular, $|B|=h-1>0$.
    
    If both $A_1$ and $A_2$ were nonempty,
    Proposition~\ref{prop:woo-component-dictionary}\ref{item:woo-two-I}
    would give two distinct elements of $\operatorname{MaxSing}(w)$,
    contradicting $|\operatorname{MaxSing}(w)|=1$.
    Hence at least one of $A_1, A_2$ is empty. 
    Since reverse-complement exchanges \(A_1\) and \(A_2\) and preserves Kazhdan--Lusztig polynomials, we may assume \(A_1=\varnothing\) by replacing \(w\), the chosen embedding, and \(\rho\) by their reverse-complements if necessary.
    We retain the same notation.
    Minimality of the amplitude also forces both
    \((a,b)\times(\alpha,\beta)\) and
    \((c,d)\times(\gamma,\delta)\) to be empty.
    Therefore,
    \[
        \bigl(w(i_1),\ldots,w(i_{h+2})\bigr)
        =
        \bigl(\alpha,\beta,\alpha-1,\ldots,\delta+1,\gamma\bigr),
    \]
    and no $w(i_j)$ lies strictly between $\gamma$ and $\delta$.
    Hence, after flattening,
    \[
        \rho(1)\cdots\rho(h+2)
        =
        (h+2)N(h+1)\cdots31.
    \]
    Since $w(d)=\delta$ occurs later, we also have $N\geq h+3$.
\end{proof}

\begin{lemma}
\label{lem:R3-parabolic-detection}
    We have
    \([q^h]P_{\mathrm{id},\rho}(q)>0\).
\end{lemma}

\begin{proof}
    Consider the standard parabolic subgroup
    \[
        W':=\langle s_1,\ldots,s_{h+2}\rangle\subseteq S_N,
    \]
    and define $y\in W'$ in one-line notation by
    \[
        y=(h+2)(h+3)(h+1)h\cdots 312(h+4)\cdots N.
    \]
    The permutation $y$ avoids $4231$, contains $3412$, and satisfies
    $h(y)=h$. Hence \([q^h]P_{\mathrm{id},y}(q)>0\) by Theorem~\ref{thm:Gaetz-Gao-min-degree}.
    
    For $g\in W'$, Lemma~\ref{lem:bruhat-tableau-criterion} gives
    \[
        g\leq\rho
        \quad\Longleftrightarrow\quad
        g(1)\leq h+2\ \text{and}\ g(h+3)\neq1
        \quad\Longleftrightarrow\quad
        g\leq y.
    \]
    For $p\geq h+3$, the first $p$ values of $g$ form
    $\{1,\ldots,p\}$, so the tableau comparisons are automatic.
    For $2\leq p\leq h+1$, the sorted prefixes of $\rho$ and $y$
    dominate every $p$-element subset of $\{1,\ldots,h+3\}$.
    The comparisons at $p=1$ and $p=h+2$ give the two conditions
    displayed above.
    Thus \([\mathrm{id},\rho]\cap W'=[\mathrm{id},y]\).
    By Theorem~\ref{thm:Billey-Braden-product},
    \(P_{\mathrm{id},y}(q)\leq_{\mathrm{coeff}}P_{\mathrm{id},\rho}(q)\).
    Hence
    \[
        [q^h]P_{\mathrm{id},\rho}(q)
        \geq
        [q^h]P_{\mathrm{id},y}(q)>0.\qedhere
    \]
\end{proof}

\begin{proof}[Proof of~\ref{case:R3}]
    By Lemmas~\ref{lem:R3-continuous-reduction} and~\ref{lem:R3-parabolic-detection}, \([q^h]P_{\mathrm{id},w}(q)>0\).
    Since $w$ contains $4231$,
    Theorem~\ref{thm:Gaetz-Gao-min-degree} gives
    \([q]P_{\mathrm{id},w}(q)>0\).
    Since $h\geq2$ and $P_{\mathrm{id},w}(1)=3$, it follows that
    \[
        P_{\mathrm{id},w}(q)=1+q+q^h.\qedhere
    \]
\end{proof}